\def\draft{n}
\documentclass[12pt]{amsart}
\usepackage[headings]{fullpage}
\usepackage{amssymb,amsmath,bbm,amsthm,graphicx}
\usepackage{epic,eepic,epsfig}
\usepackage[all]{xy}
\usepackage{url}
\usepackage[bookmarks=true,%
    colorlinks=true,%
    linkcolor=blue,%
    citecolor=blue,%
    filecolor=blue,%
    menucolor=blue,%
    urlcolor=blue,%
    breaklinks=true]{hyperref}

\newtheorem{theorem}{Theorem}[section]
\theoremstyle{definition}

\newtheorem{remark}[theorem]{Remark}

\def\printname#1{
        \if\draft y
                \smash{\makebox[0pt]{\hspace{-0.5in}
                        \raisebox{8pt}{\tt\tiny #1}}}
        \fi
}

\newlength{\standardunitlength}
\catcode`\@=11
\long\def\@makecaption#1#2{%
     \vskip 10pt

\setbox\@tempboxa\hbox{
       \small\sf{\bfcaptionfont #1. }\ignorespaces #2}%
     \ifdim \wd\@tempboxa >\captionwidth {%
         \rightskip=\@captionmargin\leftskip=\@captionmargin
         \unhbox\@tempboxa\par}%
       \else
         \hbox to\hsize{\hfil\box\@tempboxa\hfil}%
     \fi}
\font\bfcaptionfont=cmssbx10 scaled \magstephalf
\newdimen\@captionmargin\@captionmargin=2\parindent
\newdimen\captionwidth\captionwidth=\hsize
\catcode`\@=12

\def\lbl#1{\label{#1}\printname{#1}}

\def\BN{\mathbb N}
\def\BZ{\mathbb Z}
\def\BQ{\mathbbm Q}

\def\calR{\mathcal{R}}

\def\calC{\mathcal C}

\newcommand{\la}{\langle}
\newcommand{\ra}{\rangle}
\newcommand{\Ga}{\Gamma}
\newcommand{\ga}{\gamma}
\newcommand{\s}{\sigma}
\newcommand{\lla}{\langle \hspace{-2pt} \langle}
\newcommand{\rra}{\rangle \hspace{-2pt} \rangle}
\newcommand{\rot}{\mathrm{rot}}
\newcommand{\wt}{\mathrm{wt}}

\begin{document}
\title%
[A generating series for Murakami-Ohtsuki-Yamada graph evaluations]%
{A generating series for Murakami-Ohtsuki-Yamada graph evaluations}
\author{Stavros Garoufalidis}
\address{School of Mathematics \\
         Georgia Institute of Technology \\
         Atlanta, GA 30332-0160, USA \newline
         {\tt \url{http://www.math.gatech.edu/~stavros}}}
\email{stavros@math.gatech.edu}
\author{Roland van der Veen}
\address{Korteweg de Vries Institute of Mathematics \\
University of Amsterdam  
\newline 
{\tt \url{http://www.rolandvdv.nl}}}
\email{r.i.vanderveen@uva.nl}
\thanks{%
S.G. was supported in part by a National Science Foundation grant 
DMS-11-05678. R.V. was supported by the Netherlands Organization for
Scientific Research (NWO).
\bigskip\\
{\em 2010 Mathematics Subject Classification:}
  Primary 57N10. Secondary 57M25, 33F10, 39A13.\\
{\em Key words and phrases:}
  MOY graphs, colored HOMFLY polynomial, quantum topology, knots.
}

\date{March 11, 2014}

\begin{abstract}
Murakami-Ohtsuki-Yamada introduced an evaluation of certain oriented planar
trivalent graphs with colored edges. This evaluation plays a key role in the 
evaluation of the colored HOMFLY polynomial of a link in 3-space and its 
Khovanov-Rozansky categorification. Our goal is is to give a generating 
series formula for the evaluation of MOY graphs, which may be useful in
categorification, and in the study of $q$-holonomicity of the colored HOMFLY
polynomial.
\end{abstract}

\maketitle

\tableofcontents


\section{Introduction}
\lbl{sec.intro}

\subsection{The colored HOMFLY polynomial and its recursion}
\lbl{sub.intro}

The HOMFLY polynomial of a framed oriented link $L$ in 3-space is a 
powerful link invariant which takes values in the ring 
$\BZ[q^{1/2},a^{1/2},(q^{1/2}-q^{-1/2})^{-1}]$ and when specialized to
$a=q^N$, it recovers the $\mathfrak{sl}_N$ invariant of the link, 
colored by the fundamental representation. The HOMFLY polynomial has a 
colored version
$W_{L,\lambda}(q,a)$ which depends on a partition $\lambda$ for each component 
of $L$ \cite{Mo1,Mo2}. Roughly, $W_{L,\lambda}(q,a)$ is the HOMFLY of a 
universal linear combination of cables of the link $L$, where each component
colored by a partition $\lambda$ is a cabled as many times as the number
of boxes of $\lambda$ \cite{Aiston-Morton}. When suitably normalized, 
the colored HOMFLY polynomial takes values in the ring $\BZ[q^{1/2},a^{1/2}]$. 
 
In \cite{MOY}, Murakami-Ohtsuki-Yamada gave a formula for 
the colored HOMFLY polynomial in terms of evaluations of some planar,
trivalent, oriented colored graphs (in short, MOY graphs). A key 
property of a MOY graph and its evaluation is that it takes values in 
$\BN[q^{1/2}]$. The non-negativity of the coefficients of those evaluations
play an important role in categorification program developed by 
Khovanov-Rozansky \cite{KR}.

The colored HOMFLY polynomial appears in physics literature in relation to 
the large $N$ limit of $\mathrm{U}(N)$ Chern-Simons theory 
and its string dualities \cite{LMOV}. Aganagic-Vafa conjectured that
the colored HOMFLY polynomial of a knot, colored by the symmetric powers
of the fundamental representation, satisfies a recursion relation with
coefficients in $\BZ[q^{1/2},a^{1/2}]$ \cite{AV}. The operator form of such 
a recursion is a polynomial in four variables $q$, $a$, $M$ and $L$ where 
$LM=qML$ and all other variables commute. A further refinement of such an 
operator by adding a fifth variable $t$, related to the categorification of the
colored Khovanov-Rozansky Homology has been proposed by Gukov et
al \cite{Dunfield-Gukov}. Several flavors of this so-called super-polynomial
with fascinating properties have recently been conjectured in the physics 
literature. For a survey article that summarizes recent developments,
see \cite{Gukov:lectures} and \cite{AENV}.

On the mathematics side, it was observed by the first author that 
$q$-holonomicity of the colored HOMFLY polynomial (thought of as a function
of a partition with a fixed number of rows) follows from $q$-holonomicity
of the evaluations of the MOY graphs (thought of as a function of their
colors) \cite{Ga:HOMFLY}. This observation was our primary motivation to study
evaluations of MOY graphs using generating series, much in the spirit
of spin networks and their evaluations \cite{GV}. In a future publication,
we will apply our results to deduce the $q$-holonomicity of the MOY graph
evaluations.

\subsection{A generating series for the classical evaluation of 
MOY graphs}
\lbl{sub.MOY1}

\begin{figure}[htpb]
\begin{center}
\includegraphics[width=8cm]{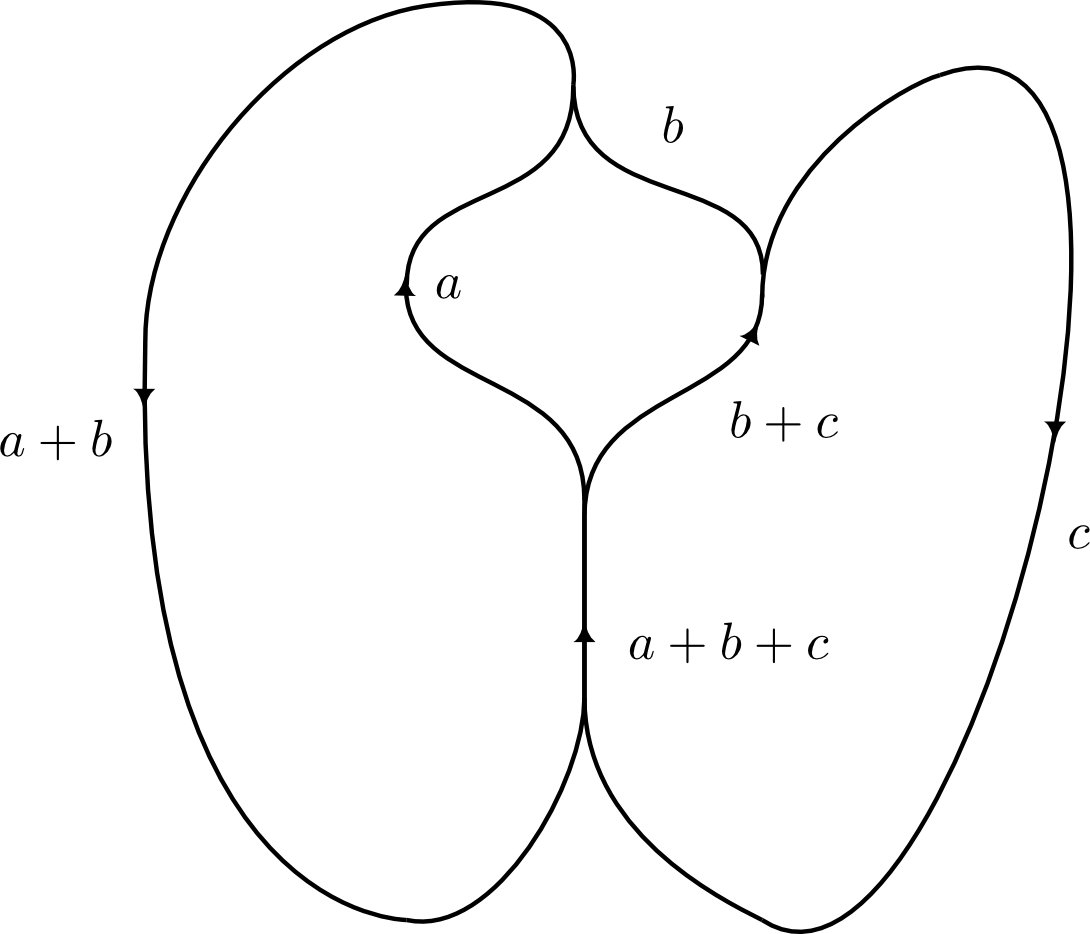}
\caption{A MOY graph and its coloring by arbitrary natural numbers 
$a$, $b$, $c$.}
\label{fig.example3}
\end{center}
\end{figure}

A MOY {\em graph} $\Ga$ is a planar trivalent graph $\Ga$ with 
oriented edges, without sinks or sources. It may contain multiple edges
and loops, as well as components with no edges. A {\em coloring} $\ga$ of
a MOY graph $\Ga$ is a {\em flow} $\ga:E(\Ga)\to \mathbb{N}$, i.e., an
assignment of a natural number to each edge such that at each vertex, the
sum of the numbers of the incoming edges equal to the 
sum of the numbers of the outgoing edges. An example of a MOY graph and
its coloring is shown in Figure \ref{fig.example3}.

For a positive natural number $N$, 
Murakami-Ohtsuki-Yamada \cite{MOY} define the evaluation 
$\la \Ga,\ga \ra_N(q) \in \BN[q^{\pm 1/2}]$. Consider the 
{\em classical evaluation} $\la \Ga,\ga \ra_N(1) \in \BN$ and its generating
series
\[
F^{\mathrm{class}}_{\Ga,N}(w) = \sum_{\ga}\la \Ga,\ga \ra_N(1) w^\ga \in \BN[w] 
\,,
\]
where $w= (w_e)_{e\in E(\Ga)}$, and all variables commute. As usual, if
$\ga: E(\Ga)\to \BN$, we denote $w^{\ga}=\prod_{\ga \in E(\Ga)} w_e^{\ga(e)}$.
The classical evaluation has been studied by Lobb-Zentner and Grant 
\cite{LZ,Grant} in relation to moduli space of $\mathfrak{sl}_N$ 
representations of the complements of MOY graphs.

A {\em cycle} of $\Ga$ is a $2$-regular subgraph of $\Ga$ such that each 
component has a consistent orientation. Let $\calC(\Ga)$ denote the
set of cycles of $\Ga$. The {\em classical cycle polynomial} is given by
$$
P^{\mathrm{class}}_\Ga(w) = \sum_{C\in \calC(\Ga)}w^C \in \BN[w] \,.
$$
Our first result identifies the generating series $F^{\mathrm{class}}_{\Ga,N}$ 
with the $N$-th power of the classical
cycle polynomial. 

\begin{theorem}
\lbl{thm.1}
We have:
\[
F^{\mathrm{class}}_{\Ga,N}(w) = \left(P^{\mathrm{class}}_\Ga(w)\right)^N \,. 
\]
\end{theorem}

\subsection{A generating series for the $\mathfrak{sl}_N$
evaluation of  MOY graphs}
\lbl{sub.MOY2}

To extend Theorem \ref{thm.1} to MOY graph evaluations, we need to
introduce the corresponding generating series and the cycle polynomial.
These are series in sets of $q$-commuting variables
$(z,Z)$ (for the generating series) and $x$ (for the cycle polynomial).

\begin{figure}[htpb]
\begin{center}
\includegraphics[width=5cm]{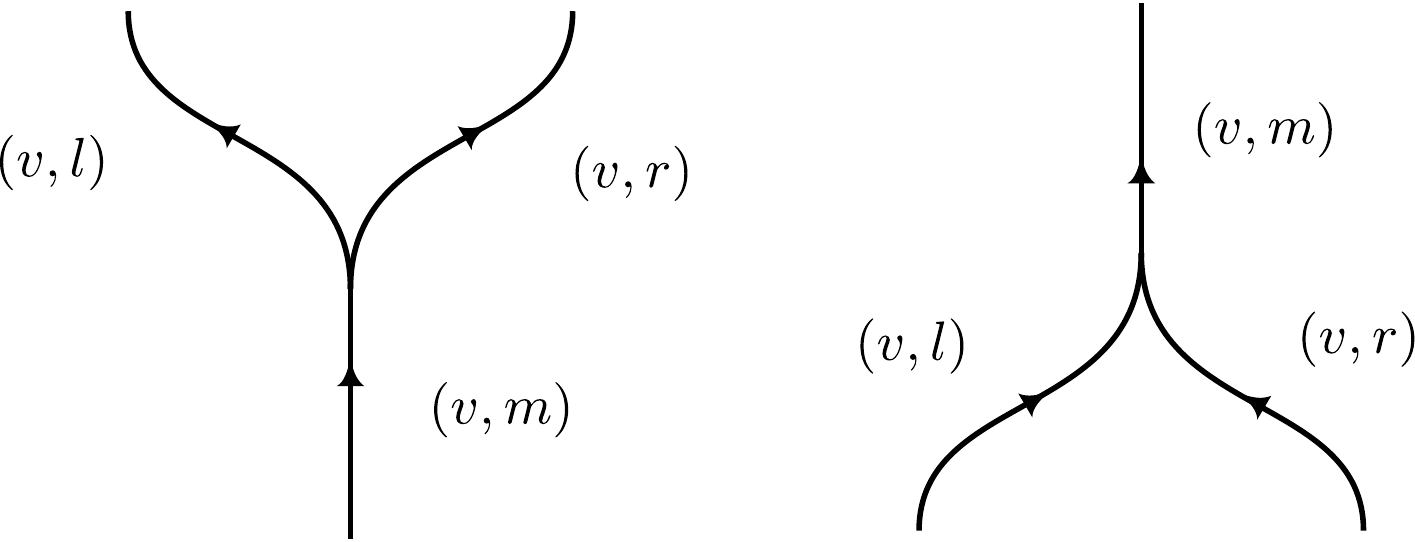}
\caption{A MOY graph and its coloring for arbitrary natural numbers 
$a$,$b$,$c$.}
\label{fig.3flags}
\end{center}
\end{figure}

To each vertex $v$ of a MOY graph $\Ga$, we denote the three adjacent 
half-edges (i.e., flags) by $(v,l)$, $(v,m)$ and $(v,r)$ with the convention
of Figure \ref{fig.3flags}. We also assign six ordered variables to $v$: 
\begin{equation}
\lbl{eq.6flags}
z_{v,l} < z_{v,m} < z_{v,r} < Z_{v,r} < Z_{v,m} < Z_{v,l}
\end{equation}
which commute except in the following instance
\begin{equation}
\lbl{eq.qv}
z_{v,r} z_{v,l} = q^{-\frac{1}{4}} z_{v,l} z_{v,r}, \qquad
Z_{v,r} Z_{v,l} = q^{-\frac{1}{4}} Z_{v,l} Z_{v,r}
\end{equation}

Fix a total ordering $<$ of the set of vertices $V(\Ga)$ of $\Ga$. Together with
\eqref{eq.6flags}, this gives a total ordering of the variables $z,Z$ where 
if $v < w$ then $z_{v,s}, Z_{v,s'} < z_{w,t}, Z_{w,t'}$ for all $s,s',t,t',
\in \{l,m,r\}$. 

Likewise, we consider a set of $q$-commuting variables 
$x=(x_C)_{C \in \calC(\Ga)}$, one for each cycle of $\Ga$. 
The commutation relations for the cycle variables are expressed in terms of 
the following intersection product (skew-symmetric form) 
on the set $\calC(\Gamma)$. 
\[
\la C,C'\ra = \frac{1}{2}\left(\#\{v\in V(\Ga)| (v,l) \in C, (v,r)\in C' \} 
- \#\{v\in V(\Ga)|(v,r) \in C, (v,l)\in C' \}\right)
\]
In terms of this product we define 
\[
x_Cx_{C'} = q^{\la C,C'\ra}x_{C'}x_{C}  \,.
\]
There is a well-defined monomial homomorphism map
\begin{equation}
\lbl{eq.mu}
\mu(x_C) =z^C Z^C \,.
\end{equation}
which satisfies $\mu(x_C x_{C'})=\mu(x_C) \mu(x_{C'})$.

The {\em generating series} of the MOY evaluations of $\Ga$ is defined by
\[
F_{\Ga,N}(q,z,Z) = \sum_{\ga}\la \Ga,\ga \ra_N(q) z^\ga Z^\ga \in 
\BZ[q^{\pm \frac{1}{2}}]\lla z,Z \rra
\,,
\]
Here the monomials $z^\ga,Z^\ga$ are understood to be in their standard ordering.

The {\em cycle polynomial} of $\Ga$ is defined in terms of the 
\emph{rotation number} $\rot(C)$ of a cycle $C$. If $C$ is connected and 
oriented counter-clockwise then $\rot(C) = 1$, if it is oriented clock-wise, 
then $\rot(C) = -1$. For a general cycle $C$, $\rot(C)$ is the sum of the 
rotation numbers of its connected components. The cycle polynomial is then
\begin{equation}
\lbl{eq.P}
P_\Ga(q,a,x) 
= \sum_{C\in \calC(\Ga)}(a^{-\frac{1}{2}}q^{\frac{1}{2}})^{\rot(C)}x_C
\in \BZ[q^{\pm \frac{1}{2}}, a^{\pm \frac{1}{2}}]\la x \ra \,.
\end{equation} 
Finally, we need to introduce a $q$-version of the $N$-th power appearing in 
Theorem \ref{thm.1}. In analogy with the $q$-Pochhammer symbol, we define
$$
(P_\Ga(q,a,x),q)_N= 
\prod_{k=0}^{N-1} P_{\Ga}(q,a,q^{k \rot(C)} x_C) \in 
\BZ[q^{\pm \frac{1}{2}}, a^{\pm \frac{1}{2}}]\la x \ra \,.
$$
We can now state our theorem.

\begin{theorem}
\lbl{thm.2}
For every MOY graph $\Gamma$ and natural number $N$ we
have:
\[
\lbl{theorem.FN}
F_{\Ga,N}(q,z,Z) = \mu \left( (P_\Ga(q,q^N,x),q)_N \right) \,.
\]
\end{theorem}
Since $\mu \left( (P_\Ga(1,1,x),1)_N \right) = (P^{\mathrm{class}}_\Ga(zZ))^N$,
Theorem \ref{thm.1} follows from Theorem \ref{thm.2} when $q=1$.
By $P^{\mathrm{class}}_\Ga(zZ)$ we mean the classical cycle polynomial where 
we set $w_e = z_{h_1}z_{h_2}Z_{h_1}Z_{h_2}$ if edge $e$ is the union of the 
half-edges $h_1, h_2$. 

\subsection{A generating series for the HOMFLY evaluation of 
MOY graphs}
\lbl{sub.MOY3}

In \cite[Lem.2.2]{Ga:HOMFLY} it was shown that given a MOY graph $(\Ga,\ga)$
there exists $\la \Ga,\ga \ra(q,a) \in \BQ(q^{1/2},a^{1/2})$ such that for every
positive natural number $N$, we have:
$$
\la \Ga,\ga \ra(q,q^N) = \la \Ga,\ga \ra_N (q) \,.
$$
Consider the generating series 
\[
F_{\Ga}(q,a,z,Z) = \sum_{\ga}\la \Ga,\ga \ra(q,a) z^\ga Z^\ga 
\]
and the rings
\[
\lbl{eq.calR}
\mathcal{R} = \BZ[[q^{\frac{1}{2}}]][a^{-\frac{1}{2}},a^{\frac{1}{2}}]
\qquad \mathcal{R}_{+} = \BZ[[q^{\frac{1}{2}}]][a^{+ \frac{1}{2}}]
\qquad \mathcal{R}_{-} = \BZ[[q^{\frac{1}{2}}]][a^{- \frac{1}{2}}]  \,.
\]
We say that a MOY graph $\Ga$ is {\em positive} if the rotation number of every 
nonempty cycle is positive. In that case, 
$F_{\Ga}(q,a,z,Z) \in \calR\lla z,Z \rra $.

\begin{theorem}
\lbl{thm.3}
Assume that $\Ga$ is positive. Then we have:
\begin{align}
\lbl{eq.FPhi}
F_{\Ga}(q,a,z,Z) \mu\left((P_{\Ga}(q,a^{-1},x),q)_\infty \right)
&= \mu\left((P_{\Ga}(q,a,x),q)_\infty\right)
\\
\lbl{eq.FPhia}
F_{\Ga}(q,q^2 a,z,Z) &= \mu\left(P_{\Ga}(q^{-1},a,x)_\infty \right) F_{\Ga}(q,a,z,Z) \mu
\left(P_{\Ga}(q,a^{-1},x)_\infty\right)
\end{align}
where
$$
(P_{\Ga}(q,a,x),q)_\infty \in \mathcal{R}_-\lla x \rra 
\qquad \mu\left((P_{\Ga}(q,a,x),q)_\infty\right) \in \mathcal{R}\lla z,Z \rra 
$$
\end{theorem}

\begin{remark}
\lbl{rem.3D}
Equation \eqref{eq.FPhi} is reminiscent to the 3D index of a 
tetrahedron introduced by Dimofte-Gaiotto-Gukov \cite{DGG2} and
further studied by \cite[Eqn.B.1]{Ga:index} and \cite{GHRS}.  
\end{remark}

\begin{remark}
\lbl{rem.HOMFLY.links}
Although non-positive MOY graphs exist (see for instance the example in
Section \ref{sub.tetrahedron}), the colored HOMFLY polynomial of a link $L$
is a linear combination, with $q$-proper hypergeometric coefficients,
of the evaluation of positive MOY graphs. Indeed, choose a braid $\beta$
whose closure is $L$ and the closure is chosen so that all strands rotate
counter-clockwise. Then, Equation on p.341 of \cite{MOY} replaces each
crossing with a linear combination of positive MOY graphs.
\end{remark}


\section{MOY graphs and their evaluation}
\lbl{sec.defs}

In this section we recall the evaluation of a MOY graph given by \cite{MOY}.

\subsection{States}

For a MOY graph $\Ga$, let $V(\Ga)$, $E(\Ga)$, $H(\Ga)$ and $\calC(\Ga)$ 
denote its set of vertices, edges, half-edges and cycles. 
For a fixed positive integer $N$, define the $N$ element set 
\[
A_N=\{-\frac{N-1}{2},\ldots,\frac{N-3}{2},\frac{N-1}{2}\}
\]

A \emph{state} is a function $\s:\calC(\Ga) \to 2^{A_N}$ 
where $2^X$ denotes the set of subsets of $X$, with the additional
requirement that if $C$ and $C'$ both contain the same edge then
$\s(C)\cap \s(C')=\emptyset$. For all states we 
require $\s(\emptyset) = \emptyset$.
A state $\s$ gives rise to functions 
$$
\s:H(\Ga)\to 2^{A_N}, \qquad
\s:E(\Ga)\to 2^{A_N}
$$
defined by $\s(h) = \bigcup_{C:h\in C} \s(C)$ and 
$\s(e) = \bigcup_{C:e\in C} \s(C)$. $\s$ induces a flow $|\s|$ on the graph 
$\Ga$ defined by $|\s|(e) = |\s(e)|$ for $e \in E(\Ga)$. 
For a cycle $C$ we define $|\s|(C) = |\s(C)|$. Finally, given a state $\s$ 
define 
\[
\mathrm{rot}(\s) = \sum_{C\in \calC(\Ga)} \rot(C)\sum_{x\in\s(C)}{x}
\] 

\subsection{Definition of the MOY evaluation}

The MOY invariant of $(\Ga,\ga)$, denoted by $\la \Ga,\ga \ra_N(q)$ 
is given by
\[
\la \Ga,\ga \ra_N(q) = \sum_{\s:|\s| = \ga} q^{\rot(\s)}\prod_{v\in V(\Ga)}\wt(v;\s)
\,. 
\]
Here define the {\em weight} by $\wt(v;\s) = q^{\frac{1}{4}(R(v;\s)-L(v;\s))}$. 
In this formula
\begin{align*}
L(v;\s) =& |\{(a,b)\in \s(v,l)\times \s(v,r)|a>b\}|\\
R(v;\s) =& |\{(a,b)\in \s(v,l)\times \s(v,r)|a<b\}|
\end{align*}
Note that in their original paper Murakami, Ohtsuki and Yamada worked with 
slightly different definitions: their concept of a state was tied to edges 
instead of cycles, but the two definitions are equivalent. Also their vertex 
weights were introduced as 
$\wt(v;\s) = q^{\frac{|\s(v,l)||\s(v,r)|-2L(v;\s)}{4}}$
which coincides with our definition above.


\section{Proofs}

In this section we present the proofs of theorems \ref{thm.2} and \ref{thm.3}.
As mentioned in the introduction theorem \ref{thm.1} follows directly from 
Theorem \ref{thm.2} by setting $q = 1$.

\subsection{Proof of Theorem \ref{thm.2}}
\lbl{sec.thm.2}

We start with the product $\mu\left((P_{\Ga}(q,q^N,x),q)_N\right)$ on the 
right hand side of the equation and show that after applying the $\mu$ map 
and ordering the variables we get the generating function
$F_{\Ga,N}(q,z,Z)$.

First we rewrite the product in a more symmetric fashion as follows:
\[
(P_{\Ga}(q,q^N,x),q)_N 
=\prod_{j\in A_N}\left(\sum_{C\in \calC(\Ga)}q^{j \rot(C)}x_C\right)
\] 

Next we need to recognize that the monomials in 
the expansion of the latter product are in bijection with the states $\s$.
Denote by $m_\s$ the monomial corresponding to state $\s$. It is defined as
\[
m_\s = \prod_{j\in A_N}\prod_{C|j\in \s(C)}x_C \,.
\]

Conversely any monomial in the expanded product looks like 
$\prod_{j\in A_N}x_{C_j}$.
This monomial corresponds to the state $\s$ defined by $\s(C) = \{j:C_j = C\}$.
Summarizing we can say that the $j$-th factor in the product corresponds to
the choice of which cycle to label by $j\in A_N$ in creating a state.

From the formula $\rot(\s)=\sum_{j\in A_N}\sum_{C:j\in \s(C)}\rot(C)$ it then
follows that
\[
(P_{\Ga}(q,q^N,x),q)_N=\sum_\s q^{\rot(\s)}m_\s \,.
\]

Our next task is to apply the monomial map $\mu$ and bring the monomials 
$\mu(m_\s)$ into the
canonical order $<$ of the $z,Z$ variables.
 We claim that the necessary $q$-commutations produce exactly coefficient
\[
\prod_{v\in V(\Ga)}\wt(v;\s) = q^{\frac{1}{4}\sum_{v}(R(v;\s)-L(v;\s))}
\] 
The $R(v;\s)$ terms come from applying $Z_{v,l}Z_{v,r} = 
q^{\frac{1}{4}}Z_{v,r}Z_{v,l}$
and the $L(v;\s)$ terms come from applying 
$z_{v,r}z_{v,l} q^{-\frac{1}{4}}z_{v,l}z_{v,r}$.
The claim now follows from the bijection between the states $\s$ and the 
monomials $m_{\s}$ because it shows that the following situations (a) and (b) 
are equivalent: 
\begin{enumerate}
\item[(a)] We have a pair of elements 
$(j_L,j_R)\in \s(C_l)\times \s(C_r)$ and a pair of cycles 
$C_l,C_r\in \calC(\Ga)$ 
such that the half-edge $(v,l)\in C_l$ and $(v,r) \in C_r$.
\item[(b)] The monomial $m_\s$ contains the factor $x_{C_l}$ in
the $j_l$-th place and $x_{C_r}$ in the $j_r$-th place. Moreover
$\mu_{C_l}$ includes $z_{v,l}Z_{v,l}$ and $\mu_{C_r}$ includes $z_{v,r}Z_{v,r}$.
\end{enumerate}
More concretely, in the graph part (a) the pair $(j_l,j_r)$ contributes $1$ to 
$L(v;\s)$ if $j_l>j_r$ and $1$ to $R(v;\s)$ otherwise.
In the monomial part (b) the case $j_l>j_r$ means that $z_r$ comes before 
$z_l$ in the product so to bring it into canonical order 
we need to commute the two and pick up a term $q^{-\frac{1}{4}}$. The case 
$j_l<j_r$ means we need to commute the upper case variables 
only and pick up a term $q^{\frac{1}{4}}$. To summarize we have now shown that
\[
\mu(m_{\s})=\prod_{v\in V(\Ga)}\wt(v;\s)z^{|\s|}Z^{|\s|}
\]
where the latter monomials are in canonical order. Therefore
 
\[
\mu\left((P_{\Ga}(q,q^N,x),q)_N\right) = \mu\left(\sum_{\s}q^{\rot(\s)}m_\s\right) 
=\sum_{\s}q^{\rot(\s)}\prod_{v\in V(\Ga)}\wt(v;\s) z^{|\s|} Z^{|\s|} \]
\[
= \sum_\ga\sum_{\s:|\s|=\ga} q^{\rot(\s)}\prod_{v\in V(\Ga)}\wt(v;\s) z^\ga Z^\ga =
F_{\Ga,N}(q,z,Z)  
\]

which concludes the proof of Theorem \ref{thm.2}.
\qed


\section{Proof of theorem \ref{thm.3}}
\lbl{sec.thm3}

For part (a) we set $a = q^N$, multiply both sides of the equality in 
Theorem \ref{thm.2} from the right by $\mu(P_{\Ga}(q,a,q^{N\rot}x),q)_\infty$ 
to obtain
\[ 
F_{\Ga}(q,a,z,Z)\mu\left((P_{\Ga}(q,a,q^{N\rot}x),q)_\infty\right) 
= \mu\left((P_{\Ga}(q,a,x),q)_\infty\right)
\]
Since
\[
P_{\Ga}(q,a,q^{N\rot}x) = 
\sum_{c\in C}(a^{\frac{1}{2}}q^{\frac{1}{2}})^{\rot(c)}x_c =  P_{\Ga}(q,a^{-1},x)
\]
we see that
\[(P_{\Ga}(q,a,q^{N\rot}x),q)_\infty = (P_{\Ga}(q,a^{-1},x),q)_\infty\]
We have now proven that for all $N$ and $a=q^N$:
\[ 
F_{\Ga}(q,a,z,Z)\mu\left((P_{\Ga}(q,a^{-1},x),q)_\infty\right) 
= \mu\left((P_{\Ga}(q,a,x),q)_\infty\right)
\]
Therefore the equality holds for all $a$.
This completes the proof of the first equality in the theorem.

For part (b) we start by writing down the statement of part (a): 
\begin{equation}
\lbl{eq.square1}
F_{\Ga}(q,a,z,Z)\mu\left((P_{\Ga}(q,a^{-1},x),q)_\infty\right) 
= \mu\left((P_{\Ga}(q,a,x),q)_\infty\right)
\end{equation}
Replacing $a$ by $q^2 a$ yields

\begin{equation}
\lbl{eq.square2}
F_{\Ga}(q,q^2 a,z,Z)\mu\left((P_{\Ga}(q,q^{-2}a^{-1},x),q)_\infty\right) =
\mu\left((P_{\Ga}(q,q^2a,x),q)_\infty\right)
\end{equation}
This can be simplified since
$P_{\Ga}(q,q^2a,x) = P_{\Ga}(q,a,q^{-\rot}x)$ and $P_{\Ga}(q,q^{-2}a^{-1},x) =
P_{\Ga}(q,a^{-1},q^{\rot}x)$
So
\begin{equation}
\lbl{eq.square3}
(P_{\Ga}(q,q^2a,x),q)_\infty = P_{\Ga}(q^{-1},a,x)(P_{\Ga}(q,a,x),q)_\infty
\end{equation}
And
\begin{equation}
\lbl{eq.square4}
P_{\Ga}(q,a^{-1},x) (P_{\Ga}(q,q^{-2}a^{-1},x),q)_\infty =
(P_{\Ga}(q,a^{-1},x),q)_\infty
\end{equation}

Substituting equation \eqref{eq.square3} into equation \eqref{eq.square2} 
we get:

\begin{equation}
\lbl{eq.square5}
F_{\Ga}(q,q^2 a,z,Z)\mu\left((P_{\Ga}(q,q^{-2}a^{-1},x),q)_\infty\right) = \mu
\left( P_{\Ga}(q^{-1},a,x)\right)\mu\left((P_{\Ga}(q,a,x),q)_\infty\right)
\end{equation}

Multiplying equation \eqref{eq.square1} from the left by 
$\mu\left(P_{\Ga}(q^{-1},a,x)\right)$ and substituting Equation 
\eqref{eq.square4} gives:

\begin{multline}
\lbl{eq.square6}
\mu\left( P_{\Ga}(q^{-1},a,x)\right) F_{\Ga}(q,a,z,Z)
\mu\left( P_{\Ga}(q,a^{-1},x)\right)
\mu\left((P_{\Ga}(q,q^{-2}a^{-1},x),q)_\infty\right) 
\\ = \mu\left( P_{\Ga}(q^{-1},a,x)\right)
\mu\left((P_{\Ga}(q,a,x),q)_\infty\right)
\end{multline}

Combining Equations \eqref{eq.square5} and \eqref{eq.square6} gives:
\begin{multline}
F_{\Ga}(q,q^2 a,z,Z)\mu\left((P_{\Ga}(q,q^{-2}a^{-1},x),q)_\infty\right) \\ = \mu
\left(P_{\Ga}(q^{-1},a,x)\right) F_{\Ga}(q,a,z,Z)\mu\left( P_{\Ga}(q,a^{-1},x)
(P_{\Ga}(q,q^{-2}a^{-1},x),q)_\infty\right)
\end{multline}
Finally multiplying by the inverse of
$\mu(P_{\Ga}(q,q^{-2}a^{-1},x),q)_\infty$ from the right we find:
\[
F_{\Ga}(q,q^2 a,z,Z) = \mu\left( P_{\Ga}(q^{-1},a,x)\right) F_{\Ga}(q,a,z,Z)
\mu\left( P_{\Ga}(q,a^{-1},x)\right)
\]
This concludes the proof of Theorem \ref{thm.3}.
\qed


\section{Examples}
\lbl{sec.examples}

\subsection{Unknot}
\lbl{sub.unknot}

For the unknot our theorems specialize to the binomial theorem and its $q$-generalization.

Let the unknot O be oriented counter-clockwise. We have no vertices and one single edge.
This means we have exactly two cycles: the empty cycle and the cycle $C$ 
that is the whole unknot.

Therefore $P_O^{class} = 1+w_C$ and so Theorem \ref{thm.1} states that
\[
\sum_{\ga}\la O,\ga\ra w^\ga= F_{O,N}^{class}(w) = (P_O^{class}(w))^N=(1+w_C)^N
\]
which is consistent with the binomial theorem and the classical evaluation
\[
\la O,\ga\ra_N(1) = {N\choose \ga}
\]
Next, according to \cite{MOY} the quantum evaluation is
\[
\la O,\ga\ra_N(q) = \left[\begin{array}{c}N\\ \ga\end{array}\right]
\]
where we are using the symmetric quantum binomial defined as
\[
\left[\begin{array}{c}n\\ k\end{array}\right] = \frac{[n]!}{[n-k]![k]!}\quad 
[n]! = \prod_{j=1}^n\frac{q^{\frac{j}{2}}-q^{-\frac{j}{2}}}{q^{\frac{1}{2}}-q^{-\frac{1}{2}}} 
\]
Theorem \ref{thm.2} now becomes a $q$-analogue of the binomial Theorem.
First note that the empty set has rotation number $0$ an the cycle $C$ rotation number $1$, 
whose cycle variable we call $x_C$. The empty set has cycle variable $1$.
Since there is only one single half-edge and no vertex we will use the 
commuting variables $z$ and $Z$ for it, so $\mu(x_C)=zZ$. 
In this notation we have
\[
P_O = 1+a^{-\frac{1}{2}}q^{\frac{1}{2}}x_C
\]
Theorem \ref{thm.2} now states that the generating series of evaluations 
\[
F_{O,N}(q,z,Z) = \sum_\ga \la O,\ga\ra(q) z^\ga Z^\ga 
= \sum_{\ga\in \mathbb{N}} \left[\begin{array}{c}N\\ \ga\end{array}\right](zZ)^\ga
\]
equals the following Pochhammer product:
\[
\mu(P_O(q,q^N,x),q)_N = \mu\prod_{k=0}^{N-1}(1+q^{\frac{1-N}{2}+k}x) 
= (-q^{\frac{1-N}{2}}zZ,q)_N
\]
The unknot is a positive so Theorem \ref{thm.3} applies, where 
$(P_O(q,a,x),q)_\infty = (-q^{\frac{1}{2}}a^{-\frac{1}{2}}x,q)_\infty$\\ 
In addition, all variables commute so we may write the theorem as
\[
F_{O}(q,a,z,Z) = \mu\left(
\frac{(P_O(q,a,x),q)_\infty}{(P_O(q,a^{-\frac{1}{2}},x),q)_\infty}\right) =
\frac{ (-q^{\frac{1}{2}}a^{-\frac{1}{2}}zZ,q)_\infty}{ (-q^{\frac{1}{2}}a^{\frac{1}{2}}zZ,q)_\infty} \,.
\]

\subsection{Tetrahedron}
\lbl{sub.tetrahedron}

\begin{figure}[htp]
\begin{center}
\includegraphics[width = 8cm]{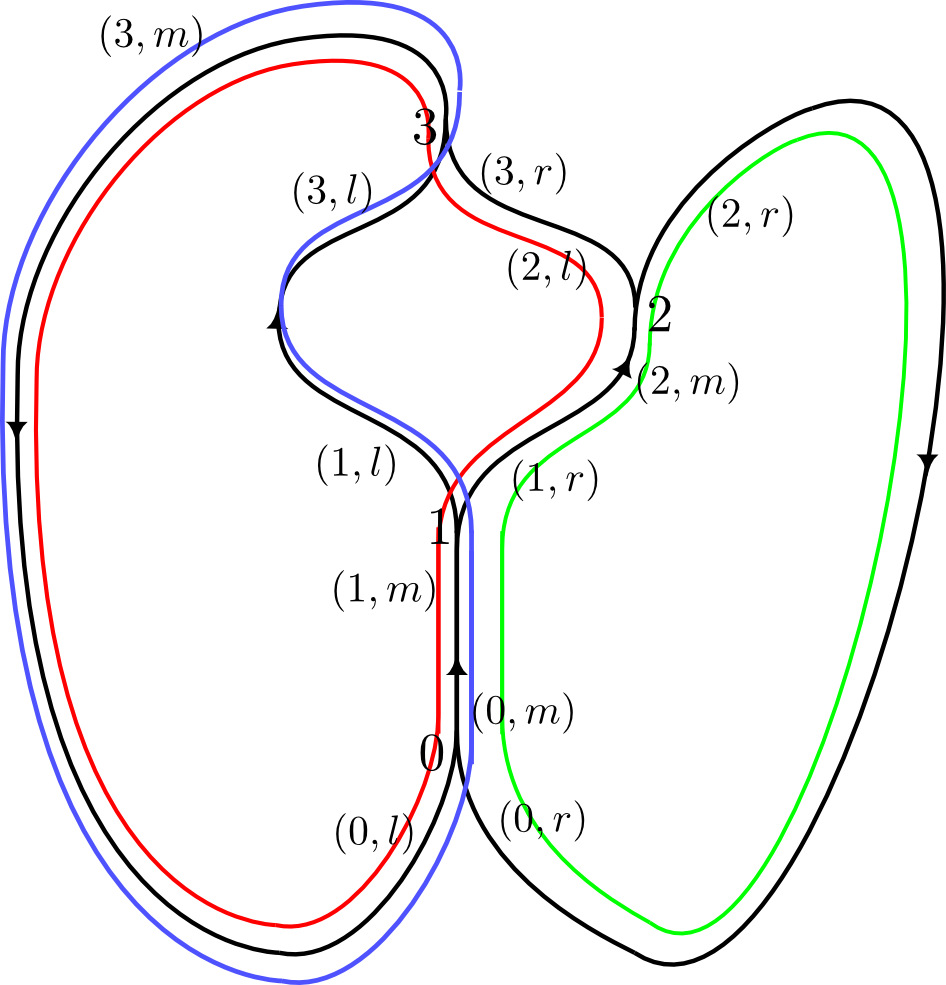}
\caption{The tetrahedron graph where we have indicated the names of the 
vertices and half edges. The three non-empty cycles are indicated in 
red, green and blue.}
\label{fig.Example2}
\end{center}
\end{figure} 

The tetrahedron graph $\Ga$ shown in Figure \ref{fig.Example2}
has four cycles, $\emptyset$ and three 
non-empty cycles $C_r$ $C_g$ and $C_b$ with classical cycle variables
$w_r$ $w_g$ and $w_b$. If $\la\Ga,a,b,c\ra_N(q)$
is the evaluation of the tetrahedron graph $\gamma$ labeled with natural numbers
$a,b,c$ as indicated in Figure \ref{fig.example3} then Theorem \ref{thm.1}
states:
\[
\sum_{a,b,c\in \mathbb{N}}\la \Ga,a,b,c\ra(1) w_r^bw_g^cw_b^a = (1+w_r+w_g+w_b)^N  
\]
in accordance to the following direct evaluation, see \cite{MOY} and the multinomial theorem.
\[
\la \Ga,a,b,c \ra_N(1) = {N \choose a,b,c,N-a-b-c}
\]
The $q$-analogue of this formula is:
\[
\la \Ga,a,b,c \ra_N(q) = \left[\begin{array}{c}N\\ a,b,c,N-a-b-c\end{array}\right]
\]

To see how Theorem \ref{thm.2} fits these numbers in a generating function, first note that
the rotation numbers of the cycles are $\rot(C_r)=\rot(C_b) = 1$ 
and $\rot(C_g)= -1$.
The corresponding cycle variables $x_r,x_g$ and $x_b$ 
can be expressed by the
$\mu$ map as follows (with the the variables in the canonical order):

\begin{align*}
\mu(x_\emptyset) &= 1\\
\mu(x_r) &= z_{0,l} z_{0,m} z_{1,m} z_{1,r} z_{2,l} z_{2,m} z_{3,r} 
z_{3,m} Z_{3,M} Z_{3,R} Z_{2,M} Z_{2,L} Z_{1,R} Z_{1,M} Z_{0,M} Z_{0,L}\\     
\mu(x_g) &= z_{0,l} z_{0,m} z_{1,l} z_{1,m} z_{3,l} z_{3,m} Z_{3,M} Z_{3,L} 
Z_{1,M} Z_{1,L} Z_{0,M} Z_{0,L}\\
\mu(x_b) &= z_{0,m} z_{0,r} z_{1,m} z_{1,r} z_{2,m} z_{2,r} Z_{2,R} Z_{2,M} 
Z_{1,R} Z_{1,M} Z_{0,R} Z_{0,M}
\end{align*}

The intersection numbers of the cycles are:
$\la x_r, x_g\ra = 1$ $\la x_r, x_b\ra = -1$ $\la x_g, x_b\ra = -1$.
From this it follows that $x_r x_g = q x_g x_r$, $x_r x_b = q^{-1} x_b x_r$ and
$x_g x_b = q^{-1} x_b x_g$. This may also be checked to follow from
applying the $\mu$ map and the commutation relations for the half-edge 
variables.

Next 
\[
P_\Ga(q,a,x) = 
1+a^{-\frac{1}{2}}q^{\frac{1}{2}}x_r+
a^{\frac{1}{2}}q^{-\frac{1}{2}}x_g+
a^{-\frac{1}{2}}q^{\frac{1}{2}}x_b
\]

The generating function is equal to $F_{\Ga,N}(q,z,Z)=$ 

\[
\sum_{a,b,c}\la \Ga, a,b,c\ra_N(q)
z_{0,l}^{a+b} z_{0,m}^{a+b+c} z_{0,r}^{c} z_{1,l}^{a} z_{1,m}^{a+b+c} z_{1,r}^{b+c} 
z_{2,l}^{b} z_{2,m}^{b+c} z_{2,r}^{c} z_{3,l}^{a} z_{3,m}^{a+b} z_{3,r}^{b}\times 
\]
\[
Z_{3,r}^{b} Z_{3,m}^{a+b} Z_{3,l}^{a} Z_{2,r}^{c} Z_{2,m}^{b+c} Z_{2,l}^{b} 
Z_{1,r}^{b+c} Z_{1,m}^{a+b+c} Z_{1,l}^{a} Z_{0,r}^{c} Z_{0,m}^{a+b+c} Z_{0,l}^{a+b} 
\]

By Theorem \ref{thm.2} this equals

\[
\mu(P_\Ga(q,q^N,x),q)_N = 
\mu\prod_{j=-\frac{N-1}{2}}^{\frac{N-1}{2}}(1+q^{j}x_r+q^{-j}x_g+q^{j}x_b)
\]

The tetrahedron graph $\Ga$ in our example is not positive, $\rot(C_g)=-1$, so Theorem \ref{thm.3} does not apply. 
However the green cycle may be turned into a positive cycle by moving its 
ends to the left.

For several reasons this example may still be too simple in that there
are no relations between the cycles. For more complicated graphs, monomials 
can usually be written as a product of cycles in multiple ways. Also a simple closed form
evaluation in terms of $q$-binomials is generally not to be expected.


\subsection*{Acknowledgment}
The first author wishes to thank Christoph Koutschan for enlightening
conversations.

\bibliographystyle{hamsalpha}
\bibliography{biblio}
\end{document}